\documentclass[10pt, a4paper]{article}

\usepackage{amsmath,  amsfonts}

\usepackage[T1]{fontenc}

\newtheorem{Theorem}{\quad Theorem}[section]

\newcommand{\be} {\begin{equation}}
\newcommand{\ee} {\end{equation}}

\title {A compactness result for solutions to an equation with boundary singularity.}

\date{}

 \author{Samy Skander Bahoura\footnote {e-mails: samybahoura@yahoo.fr, samybahoura@gmail.com} \\ 
 {\small Equipe d'Analyse Complexe et G\'eom\'etrie.}\\  
  {\small Universit\'e Pierre et Marie Curie, 75005 Paris, France.}} 

\begin{document}

\maketitle
 
\begin{abstract}

We give a blow-up behavior for solutions to a problem with boundary singularity and Dirichlet condition. An application, we derive a compactness result for solutions to this Problem with singularity and Lipschitz condition. 

\end{abstract}

{ \small  Keywords: blow-up, boundary, logarithmic singularity, a priori estimate, analytic domain, Lipschitz condition.}

{\bf \small MSC: 35J60, 35J05, 35B44, 35B45}

\section{Introduction and Main Results} 

We set $ \Delta = \partial_{11} + \partial_{22} $  on an analytic domain $ \Omega \subset {\mathbb R}^2 $.

\bigskip

We consider the following equation:

$$ (P)   \left \{ \begin {split} 
      -\Delta u & = -\log \dfrac{|x|}{2d}V e^{u} \,\, &\text{in} \,\, & \Omega  \subset {\mathbb R}^2, \\
                  u & = 0  \,\,             & \text{in} \,\,    &\partial \Omega.              
\end {split}\right.
$$

Here:

$$ d=diam(\Omega), \,\, 0 \in \partial \Omega, $$

and,
 
$$ 0\leq V \leq b <+\infty,\,\,  -\log \dfrac{|x|}{2d}e^u \in L^1(\Omega),\,\, u\in W^{1,1}_0(\Omega). $$

The previous equation was studied by many authors, with or without  the boundary condition, also for Riemannian surfaces,  see [1-15],  we can find some existence and compactness results.

Among other results, we  can see in [7] the following important Theorem,

\smallskip

{\bf Theorem.}{\it (Brezis-Merle [7])}.{\it If $ (u_i)_i $ and $ (V_i)_i $ are two sequences of functions relatively to the problem $ (P) $ with, $ 0 < a \leq V_i \leq b < + \infty $, then, for all compact set $ K $ of $ \Omega $,

$$ \sup_K u_i \leq c = c(a, b, K, \Omega). $$}

If we assume $ V $ with more regularity, we can have another type of estimates, a $ \sup + \inf $ type inequalities. It was proved by Shafrir see [15], that, if $ (u_i)_i, (V_i)_i $ are two sequences of functions solutions of the previous equation without assumption on the boundary and, $ 0 < a \leq V_i \leq b < + \infty $, then we have the following interior estimate:

$$ C\left (\dfrac{a}{b} \right ) \sup_K u_i + \inf_{\Omega} u_i \leq c=c(a, b, K, \Omega). $$

\bigskip

Now, if we suppose $ (V_i)_i $ uniformly Lipschitzian with $ A $ the
Lipschitz constant, then, $ C(a/b)=1 $ and $ c=c(a, b, A, K, \Omega)
$, see [5]. 
\smallskip

Here we give the behavior of the blow-up points on the boundary and a proof of a compactness result of solutions of a Brezis-Merle type Problem with Lipschitz condition.

\smallskip

Here, we write an extenstion of Brezis-Merle Problem (see [7]) is:

\smallskip

{\bf Problem}. Suppose that $ V_i \to  V $ in $ C^0( \bar \Omega ) $, with, $ 0 \leq V_i $. Also, we consider a sequence of solutions $ (u_i) $ of $ (P) $ relatively to $ (V_i) $ such that,

$$ \int_{\Omega} -\log \dfrac{|x|}{2d} e^{u_i} dx \leq C,  $$

is it possible to have:

$$ ||u_i||_{L^{\infty}}\leq C ? $$

Here, we give a caracterization of the behavior of the blow-up points on the boundary and also a proof of the compactness theorem when the prescribed curvature are uniformly Lipschitzian. For the behavior of the blow-up points on the boundary, the following condition is enough,

$$ 0 \leq  V_i \leq b, $$

The condition $ V_i \to  V $ in $ C^0(\bar \Omega) $ is not necessary.

\bigskip

But for the proof  for the Brezis- Merle type problem we assume that:

$$ ||\nabla V_i||_{L^{\infty}}\leq  A. $$

We have the following caracterization of the behavior of the blow-up points on the boundary.

\begin{Theorem}  Assume that $ \max_{\Omega} u_i \to +\infty $, Where $ (u_i) $ are solutions of the probleme $ (P) $ with:
 
 $$ 0 \leq V_i \leq b,\,\,\, {\rm and } \,\,\, \int_{\Omega}  -\log \dfrac{|x|}{2d}  e^{u_i} dx \leq C, \,\,\, \forall \,\, i, $$
 then;  after passing to a subsequence, there is a finction $ u $,  there is a number $ N \in {\mathbb N} $ and there are $ N  $ points $ x_1, x_2, \ldots, x_N \in  \partial \Omega $, such that, 

$$ \partial_{\nu} u_i  \to \partial_{\nu} u +\sum_{j=1}^N \alpha_j \delta_{x_j}, \,\,\, \alpha_1 \geq 4\pi,\alpha_j \geq 4\pi \,\, {\rm weakly\,\, in \, the \, sens \, of \, measure } \,\, L^1(\partial\Omega). $$

and,

$$ u_i \to u \,\,\, {\rm in }\,\,\, C^1_{loc}(\bar \Omega-\{x_1,\ldots, x_N \}). $$

\end{Theorem} 

 In the following theorem, we have a proof for the global a priori estimate which concern the problem $ (P) $.

\bigskip

\begin{Theorem}Assume that $ (u_i) $ are solutions of $ (P) $ relatively to $ (V_i) $ with the following conditions:

$$ d=diam (\Omega),\,\, 0 \in \partial \Omega, $$

$$  0 \leq V_i \leq b, \,\, ||\nabla V_i||_{L^{\infty}} \leq A,\,\, {\rm and } \,\,\, \int_{\Omega} -\log \dfrac{|x|}{2d}e^{u_i} \leq C, $$

We have,

$$  || u_i||_{L^{\infty}} \leq c(b, A, C, \Omega), $$

\end{Theorem}

\section{Proof of the theorems} 

\bigskip

\underbar {\it Proof of theorem 1.1:} 

\bigskip

We have,

$$ u_i\in W^{1,1}_0(\Omega). $$

By [7], $ e^{u_i} \in L^k,\,\, \forall k >1 $ and by the elliptic estimates:

$$ u_i\in W^{2,k}(\Omega) \cap C^{1,\epsilon}(\bar \Omega),\,\, \epsilon >0. $$

We have, 

$$ \int_{\partial \Omega} \partial_{\nu} u_i d\sigma \leq C, $$

Without loss of generality, we can assume that $ \partial_{\nu} u_i \geq 0 $. Thus, (using the weak convergence in the space of Radon measures), we have the existence of a positive Radon measure $ \mu $ such that,

$$ \int_{\partial \Omega} \partial_{\nu} u_i \phi  d\sigma \to \mu(\phi), \,\,\, \forall \,\,\, \phi \in C^0(\partial \Omega). $$

We take an $ x_0 \in \partial \Omega $ such that, $ \mu({x_0}) < 4\pi $. Without loss of generality, we can assume that the following curve, $ B(x_0, \epsilon) \cap \partial \Omega := I_{\epsilon} $ is an interval.(In this case, it is more simple to construct the following test function $ \eta_{\epsilon} $). We choose a function $ \eta_{\epsilon} $ such that,

$$ \begin{cases}
    
\eta_{\epsilon} \equiv 1,\,\,\,  {\rm on } \,\,\,  I_{\epsilon}, \,\,\, 0 < \epsilon < \delta/2,\\

\eta_{\epsilon} \equiv 0,\,\,\, {\rm outside} \,\,\, I_{2\epsilon }, \\

0 \leq \eta_{\epsilon} \leq 1, \\

||\nabla \eta_{\epsilon}||_{L^{\infty}(I_{2\epsilon})} \leq \dfrac{C_0(\Omega, x_0)}{\epsilon}.

\end{cases} $$

We take a $\tilde \eta_{\epsilon} $ such that,

$$  \left \{ \begin {split} 
       -\Delta \tilde \eta_{\epsilon} &= 0 && \text{in} \,\,\Omega \\
            \tilde\eta_{\epsilon} &=  \eta_{\epsilon} && \text{on} \,\, \partial \Omega.              
\end {split}\right.
$$

{\bf Remark:} We use the following steps in the construction of $ \tilde \eta_{\epsilon} $:

We take a cutoff function $ \eta_{0} $ in $ B(0, 2) $ or $ B(x_0, 2) $:

We set $ \eta_{\epsilon}(x)= \eta_0(|x-x_0|/\epsilon) $ in the case of the unit disk it is sufficient.

\smallskip

We use a chart $ [f, B_1(0)] $ with $ f(0)=x_0 $.

\smallskip

We can take: $ \mu_{\epsilon}(x)= \eta_0(x/\epsilon) $ and $ \eta_{\epsilon}(y)= \mu_{\epsilon}(f^{-1}(y)) $, we extend it by $ 0 $ outside $ f(B_1(0)) $.  We have $ f(B_1(0)) = D_1(x_0) $, $ f (B_{\epsilon}(0))= D_{\epsilon}(x_0) $ and $ f(B_{\epsilon}^+)= D_{\epsilon}^+(x_0) $ with $ f $ and $ f^{-1} $ smooth diffeomorphism.

$$ \begin{cases}
    
\eta_{\epsilon} \equiv 1,\,\,\,  {\rm on \, a \, the \, connected \, set } \,\,\,  J_{\epsilon} =f(I_{\epsilon}), \,\,\, 0 < \epsilon < \delta/2,\\

\eta_{\epsilon} \equiv 0,\,\,\, {\rm outside} \,\,\, J'_{\epsilon}=f(I_{2\epsilon }), \\

0 \leq \eta_{\epsilon} \leq 1, \\

||\nabla \eta_{\epsilon}||_{L^{\infty}(J'_{\epsilon})} \leq \dfrac{C_0(\Omega, x_0)}{\epsilon}.

\end{cases} $$

And, $ H_1(J'_{\epsilon}) \leq C_1 H_1(I_{2\epsilon}) = C_1 4\epsilon $, since $ f $ is Lipschitz. Here $ H_1 $ is the Hausdorff measure.

 We solve the Dirichlet Problem:

\begin{displaymath}  \left \{ \begin {split} 
      \Delta \bar \eta_{\epsilon}  & = \Delta \eta_{\epsilon}              \,\, &&\text{in} \!\!&&\Omega \subset {\mathbb R}^2, \\
                  \bar \eta_{\epsilon} & = 0   \,\,             && \text{in} \!\!&&\partial \Omega.               
\end {split}\right.
\end{displaymath}

and finaly we set $ \tilde \eta_{\epsilon} =-\bar \eta_{\epsilon} + \eta_{\epsilon} $. Also, by the maximum principle and the elliptic estimates we have :

$$ ||\nabla \tilde \eta_{\epsilon}||_{L^{\infty}} \leq C(|| \eta_{\epsilon}||_{L^{\infty}} +||\nabla \eta_{\epsilon}||_{L^{\infty}} + ||\Delta \eta_{\epsilon}||_{L^{\infty}}) \leq \dfrac{C_1}{\epsilon^2}, $$

with $ C_1 $ depends on $ \Omega $.

We use the following estimate, see [8],

$$ ||\nabla u_i||_{L^q} \leq C_q, \,\,\forall \,\, i\,\, {\rm and  }  \,\, 1< q < 2. $$

We deduce from the last estimate that, $ (u_i) $ converge weakly in $ W_0^{1, q}(\Omega) $, almost everywhere to a function $ u \geq 0 $ and $ \int_{\Omega} -\log \dfrac{|x|}{2d}e^u < + \infty $ (by Fatou lemma). Also, $ V_i $ weakly converge to a nonnegative function $ V $ in $ L^{\infty} $. The function $ u $ is in $ W_0^{1, q}(\Omega) $ solution of :

$$  \left \{ \begin {split} 
       -\Delta u &= -\log \dfrac{|x|}{2d}Ve^{u} \in L^1(\Omega) && {\rm in } \,\, \Omega \\
                   u&= 0 && \text{on} \,\,\partial \Omega,              
\end {split}\right.
$$

As in  the corollary 1 of Brezis-Merle result, see [7],   we have $ e^{k u }\in L^1(\Omega), k >1 $. By the elliptic estimates, we have $ u \in C^1(\bar \Omega) $.

\bigskip

We can write,

 \be -\Delta ((u_i-u) \tilde \eta_{\epsilon})= -\log\dfrac{|x|}{2d}(V_i e^{u_i} -Ve^u)\tilde \eta_{\epsilon} +2(\nabla (u_i- u) \cdot \nabla \tilde \eta_{\epsilon}) .\label{(1)}\ee

We use the interior esimate of Brezis-Merle, see [7],

\bigskip

\underbar {\it Step 1:} Estimate of the integral of the first term of the right hand side of $ (\ref{(1)}) $.

\bigskip

We use the Green formula between $ \tilde \eta_{\epsilon} $ and $ u $, we obtain,

 \be \int_{\Omega} -\log \dfrac{|x|}{2d}Ve^u \tilde \eta_{\epsilon} dx =\int_{\partial \Omega} \partial_{\nu} u \eta_{\epsilon} \leq 4\epsilon ||\partial_{\nu}u||_{L^{\infty}}= C \epsilon \label{(2)}\ee

We have,

$$  \left \{ \begin {split} 
       -\Delta u_i &= -\log \dfrac{|x|}{2d}V_ie^{u_i} && {\rm in } \,\, \Omega \\
                   u_i&= 0 && \text{on} \,\,\partial \Omega,              
\end {split}\right.
$$

We use the Green formula between $ u_i $ and $ \tilde \eta_{\epsilon} $ to have:

\be \int_{\Omega} -\log \dfrac{|x|}{2d}V_i e^{u_i} \tilde \eta_{\epsilon} dx = \int_{\partial \Omega} \partial_{\nu} u_i \eta_{\epsilon} d\sigma \to \mu(\eta_{\epsilon}) \leq \mu(I_{2\epsilon}) \leq 4  \pi - \epsilon_0, \,\,\, \epsilon_0 >0 \label{(3)}\ee

From $(\ref{(2)}) $ and $ (\ref{(3)}) $ we have for all $ \epsilon >0 $ there is $ i_0 =i_0(\epsilon) $ such that, for $ i \geq i_0 $,

\be \int_{\Omega} -\log \dfrac{|x|}{2d}|(V_ie^{u_i}-Ve^u) \tilde \eta_{\epsilon}| dx \leq 4 \pi -\epsilon_0 + C \epsilon \label{(4)} \ee

\bigskip

\underbar { Step 2:} Estimate of integral of the second term of the right hand side of $ (\ref{(1)}) $.

\bigskip
Let $ \Sigma_{\epsilon} = \{ x \in \Omega, d(x, \partial \Omega) = \epsilon^3 \} $ and $ \Omega_{\epsilon^3} = \{ x \in \Omega, d(x, \partial \Omega) \geq \epsilon^3 \} $, $ \epsilon > 0 $. Then, for $ \epsilon $ small enough, $ \Sigma_{\epsilon} $ is hypersurface.

\bigskip

The measure of $ \Omega-\Omega_{\epsilon^3} $ is $ k_2\epsilon^3 \leq \mu_L (\Omega-\Omega_{\epsilon^3}) \leq k_1 \epsilon^3 $.

\bigskip

{\bf Remark}: for the unit ball $ {\bar B(0,1)} $, our new manifold is $ {\bar B(0, 1-\epsilon^3)} $.

\bigskip

(Proof of this fact; let's consider $ d(x, \partial \Omega) = d(x, z_0), z_0 \in \partial \Omega $, this imply that $ (d(x,z_0))^2 \leq (d(x, z))^2 $ for all $ z \in \partial \Omega $ which it is equivalent to $ (z-z_0) \cdot (2x-z-z_0) \leq 0 $ for all $ z \in \partial \Omega $, let's consider a chart around $ z_0 $ and $ \gamma (t) $ a curve in $ \partial \Omega $, we have;

$ (\gamma (t)-\gamma(t_0) \cdot (2x-\gamma(t)-\gamma(t_0)) \leq 0 $ if we divide by $ (t-t_0) $ (with the sign and tend $ t $ to $ t_0 $), we have $ \gamma'(t_0) \cdot (x-\gamma(t_0)) = 0 $, this imply that $ x= z_0-s \nu_0 $ where $ \nu_0 $ is the outward normal of $ \partial \Omega $ at $ z_0 $))

With this fact, we can say that $ S= \{ x, d(x, \partial \Omega) \leq \epsilon \}= \{ x= z_0- s \nu_{z_0}, z_0 \in \partial \Omega, \,\, -\epsilon \leq s \leq \epsilon \} $. It  is sufficient to work on  $ \partial \Omega $. Let's consider a charts $ (z, D=B(z, 4 \epsilon_z), \gamma_z) $ with $ z \in \partial \Omega $ such that $ \cup_z B(z, \epsilon_z) $ is  cover of $ \partial \Omega $ .  One can extract a finite cover $ (B(z_k, \epsilon_k)), k =1, ..., m $, by the area formula the measure of $ S \cap B(z_k, \epsilon_k) $ is less than a $ k\epsilon $ (a $ \epsilon $-rectangle).  For the reverse inequality, it is sufficient to consider one chart around one point of the boundary).

\bigskip

We write,

$$ \int_{\Omega} |(\nabla ( u_i -u) \cdot \nabla \tilde \eta_{\epsilon} ) | dx =
\int_{\Omega_{\epsilon^3}} |(\nabla (u_i -u) \cdot \nabla \tilde \eta_{\epsilon}) | dx + $$

\be + \int_{\Omega - \Omega_{\epsilon^3}} (\nabla (u_i-u) \cdot \nabla \tilde \eta_{\epsilon})| dx. \label{(5)} \ee

\bigskip

\underbar {\it Step 2.1:} Estimate of $ \int_{\Omega - \Omega_{\epsilon^2}} |(\nabla (u_i-u) \cdot \nabla \tilde \eta_{\epsilon})| dx $.

\bigskip
First, we know from the elliptic estimates that  $ ||\nabla \tilde \eta_{\epsilon}||_{L^{\infty}} \leq C_1 /\epsilon^2 $, $ C_1 $ depends on $ \Omega $

\bigskip

We know that $ (|\nabla u_i|)_i $ is bounded in $ L^q, 1< q < 2 $, we can extract  from this sequence a subsequence which converge weakly to $ h \in L^q $. But, we know that we have locally the uniform convergence to $ |\nabla u| $ (by Brezis-Merle theorem), then, $ h= |\nabla u| $ a.e. Let $ q' $ be the conjugate of $ q $.

\bigskip

We have, $  \forall f \in L^{q'}(\Omega)$

$$ \int_{\Omega} |\nabla u_i| f dx \to \int_{\Omega} |\nabla u| f dx $$

\bigskip

If we take $ f= 1_{\Omega-\Omega_{\epsilon^3}} $, we have:

$$ {\rm for } \,\, \epsilon >0 \,\, \exists \,\, i_1 = i_1(\epsilon) \in {\mathbb N}, \,\,\, i \geq  i_1,  \,\, \int_{\Omega-\Omega_{\epsilon^3}} |\nabla u_i| \leq \int_{\Omega-\Omega_{\epsilon^2}} |\nabla u| + \epsilon^2. $$

Then, for $ i \geq i_1(\epsilon) $,

$$ \int_{\Omega-\Omega_{\epsilon^3}} |\nabla u_i| \leq mes(\Omega-\Omega_{\epsilon^3}) ||\nabla u||_{L^{\infty}} + \epsilon^2 = C \epsilon. $$

Thus, we obtain,

\be \int_{\Omega - \Omega_{\epsilon^3}} |(\nabla (u_i-u) \cdot \nabla \tilde \eta_{\epsilon})| dx \leq  \epsilon C_1(2 k_1 ||\nabla u||_{L^{\infty}} + 1) \label{(6)}\ee

The constant $ C_1 $ does  not depend on $ \epsilon $ but on $ \Omega $.

\bigskip

\underbar {\it Step 2.2:} Estimate of $ \int_{\Omega_{\epsilon^2}} |(\nabla (u_i-u) \cdot \nabla \tilde \eta_{\epsilon})| dx $.

\bigskip

We know that, $ \Omega_{\epsilon} \subset \subset \Omega $, and ( because of Brezis-Merle's interior estimates) $ u_i \to u $ in $ C^1(\Omega_{\epsilon^3}) $. We have,

$$ ||\nabla (u_i-u)||_{L^{\infty}(\Omega_{\epsilon^3})} \leq \epsilon^2,\, {\rm for } \,\, i \geq i_3 = i_3(\epsilon). $$

We write,
 
$$ \int_{\Omega_{\epsilon_3}} |(\nabla (u_i-u) \cdot \nabla \tilde \eta_{\epsilon} )| dx \leq ||\nabla (u_i-u)||_{L^{\infty}(\Omega_{\epsilon^3})} ||\nabla \tilde \eta_{ \epsilon}||_{L^{\infty}} \leq C_1 \epsilon \,\, {\rm for } \,\, i \geq i_3, $$

For $ \epsilon >0 $, we have for $ i \in {\mathbb N} $, $ i \geq \max \{i_1, i_2, i_3 \} $,

\be \int_{\Omega} |(\nabla (u_i-u) \cdot \nabla \tilde \eta_{\epsilon})| dx \leq \epsilon C_1(2 k_1 ||\nabla u||_{L^{\infty}} + 2) \label{(7)}\ee

From $ (\ref{(4)}) $ and $(\ref{(7)}) $, we have, for $ \epsilon >0 $, there is $ i_3= i_3(\epsilon) \in {\mathbb N}, i_3 = \max \{ i_0, i_1, i_2 \} $ such that,

\be \int_{\Omega} |\Delta [(u_i-u)\tilde \eta_{\epsilon}]|dx \leq 4 \pi-\epsilon_0+  \epsilon 2 C_1(2 k_1 ||\nabla u||_{L^{\infty}} + 2 + C) \label{(8)}\ee

We choose $ \epsilon >0 $ small enough to have a good estimate of  $ (1) $.

\bigskip

Indeed, we have:

$$ \left \{ \begin {split}
 -\Delta [(u_i-u) \tilde \eta_{\epsilon}] &= g_{i,\epsilon} && {\rm in } \,\,\, \Omega,\\
 (u_i-u) \tilde \eta_{\epsilon} &= 0 && {\rm on } \,\,\, \partial \Omega.
 \end {split}  \right . $$

with $ ||g_{i, \epsilon} ||_{L^1(\Omega)} \leq 4 \pi -\epsilon_0. $

\bigskip

We can use Theorem 1 of [7] to conclude that there is $ q \geq\tilde q >1 $ such that:

$$ \int_{V_{\epsilon}(x_0)} e^{\tilde q|u_i-u|} dx \leq \int_{\Omega} e^{q|u_i -u| \tilde \eta_{\epsilon}} dx \leq C(\epsilon,\Omega). $$
 
where, $ V_{\epsilon}(x_0) $ is a neighberhood of $ x_0 $ in $ \bar \Omega $.

\bigskip

Thus, for each $ x_0 \in \partial \Omega - \{ \bar x_1,\ldots, \bar x_m \} $ there is $ \epsilon_{x_0} >0, q_{x_0} > 1 $ such that:

$$ \int_{B(x_0, \epsilon_{x_0})} e^{q_{x_0} u_i} dx \leq C, \,\,\, \forall \,\,\, i. $$

Now, we consider a cutoff function $ \eta \in C^{\infty}({\mathbb R}^2) $ such that:

$$ \eta \equiv 1 \,\,\, {\rm on } \,\,\, B(x_0, \epsilon_{x_0}/2) \,\,\, {\rm and } \,\,\, \eta \equiv 0 \,\,\, {\rm on } \,\,\, {\mathbb R}^2 -B(x_0, 2\epsilon_{x_0}/3). $$

We write,

$$ -\Delta (u_i \eta) = -\log \dfrac{|x|}{2d}V_i e^{u_i} \eta - 2 <\nabla u_i|\nabla \eta> - u_i \Delta \eta. $$

By the elliptic estimates, $ (u_i \eta)_i $ is uniformly bounded in $ W^{2, q_1}(\Omega) $ and also, in $ C^1(\bar \Omega) $.

\bigskip

Finaly, we have, for some $ \epsilon > 0 $ small enough,

\be || u_i||_{C^{1,\theta}[B(x_0, \epsilon)]} \leq c_3 \,\,\, \forall \,\,\, i. \label{(9)}\ee

We have proved that, there is a finite number of points $ \bar x_1, \ldots, \bar x_m $ such that the squence $ (u_i)_i  $ is locally uniformly bounded in $ \bar \Omega - \{ \bar x_1, \ldots , \bar x_m \} $.

\smallskip

And, finaly, we have:

\be \partial_{\nu} u_i  \to \partial_{\nu} u +\sum_{j=1}^N \alpha_j \delta_{x_j}, \,\,\alpha_j \geq 4\pi \,\, {\rm weakly\,\, in \, the \, sens \, of \, measure } \,\, L^1(\partial\Omega). \ee

\underbar {\it Proof of theorem 1.2:} 

\bigskip
Without loss of generality, we can assume that $ 0 $ is a blow-up point. Since the boundary is an analytic curve $ \gamma(t) $, there is a neighborhood of  $ 0=x_1 $ such that the curve $ \gamma $ can be extend to a holomorphic map such that $ \gamma'(0) \not = 0 $ (series) and by the inverse mapping one can assume that this map is univalent around $ 0 $. In the case when the boundary is a simple Jordan curve the domain is simply connected. In the case that the domains has a finite number of holes it is conformally equivalent to a disk with a finite number of disks removed. Here we consider a general domain. Without loss of generality one can assume that $ \gamma (B_1^+) \subset \Omega $ and also $ \gamma (B_1^-) \subset (\bar \Omega)^c $ and $ \gamma (-1,1) \subset \partial \Omega $ and $ \gamma $ is univalent. This means that $ (B_1, \gamma) $ is a local chart around $ 0 $ for $ \Omega $ and $ \gamma $ univalent. (This fact holds if we assume that we have an analytic domain, in the sense of Hofmann see [11], (below a graph of an analytic function), we have necessary the condition $ \partial \bar \Omega = \partial \Omega $ and the graph is analytic, in this case $ \gamma (t)= (t, \phi(t)) $ with $ \phi $ real analytic and an example of this fact is the unit disk  around the point $ (0,1) $ for example).

By this conformal transformation, we can assume that $ \Omega =B_1^+ $, the half ball, and $ \partial^+ B_1^+ $ is the exterior part, a part which not contain $ 0 $ and on which  $ u_i $ converge in the $ C^1 $ norm to $ u $. Let us consider $ B_{\epsilon}^+ $, the half ball with radius $ \epsilon >0 $. Also, one can consider a $ C^1 $ domain (a rectangle between two half disks) and by charts its image is a $ C^1 $ domain).

We know that:

$$ u_i \in W^{2,k}(\Omega)\cap C^{1, \epsilon}(\bar \Omega). $$ 

Thus we can use integrations by parts (Gauss-Green-Riemann-Stokes formula). The second Pohozaev identity applied around the blow-up $ 0=x_1 $ gives :

\be  \int_{B_{\epsilon}^+} -\Delta u_i (x \cdot \nabla u_i) dx = -\int_{\partial^+ B_{\epsilon}^+}  g(\nabla u_i)d\sigma, \label{(10)}\ee

with,

$$ g(\nabla u_i)= (\nu \cdot \nabla u_i)(x \cdot\nabla u_i)- x \cdot \nu \dfrac{|\nabla u_i|^2}{2}. $$

After integration by parts, we obtain:

$$  \int_{B_{\epsilon}^+} -\log \dfrac{|x|}{2d}V_ie^{u_i}(1+o(1)) dx + \int_{B_{\epsilon}^+}  (x \cdot \nabla V_i)(-\log \dfrac{|x|}{2d}) V_ie^{u_i} dx + \int_{\partial B_{\epsilon}^+} -\log \dfrac{|x|}{2d}V_ie^{u_i}(x \cdot \nu) =$$

\be \int_{\partial^+ B_{\epsilon}^+}  g(\nabla u_i)d\sigma +o(1), \label{(10)}\ee

$$ \int_{B_{\epsilon}^+} -\log \dfrac{|x|}{2d} Ve^{u}(1+o(1)) dx + \int_{B_{\epsilon}^+}  (x \cdot \nabla V)(-\log \dfrac{|x|}{2d}) Ve^{u} dx +\int_{\partial B_{\epsilon}^+} -\log \dfrac{|x|}{2d}Ve^{u}(x \cdot \nu) =$$
 
\be = \int_{\partial^+ B_{\epsilon}^+}  g(\nabla u)d\sigma +o(1), \label{(10)}\ee

Thus, 

$$ \int_{B_{\epsilon}^+} -\log \dfrac{|x|}{2d} V_ie^{u_i} dx-\int_{B_{\epsilon}^+} -\log \dfrac{|x|}{2d}Ve^{u} dx + $$

$$ + \int_{B_{\epsilon}^+}  (x \cdot \nabla V_i)(-\log \dfrac{|x|}{2d})V_ie^{u_i} dx - \int_{B_{\epsilon}^+}  (x \cdot \nabla V)(-\log \dfrac{|x|}{2d})Ve^{u} dx = $$

$$  = \int_{\partial^+ B_{\epsilon}^+}  g(\nabla u_i)-g(\nabla u) d\sigma +o(1)= o(1), $$

First, we tend $ i  $ to infinity after $ \epsilon  $ to 0, we obtain:

\be \lim_{ \epsilon \to 0} \lim_{ i\to + \infty}\int_{B_{\epsilon}^+} -\log \dfrac{|x|}{2d}V_ie^{u_i} dx = 0, \label{(12)}\ee

But,

$$ \int_{\gamma (B_{\epsilon}^+)} -\log \dfrac{|x|}{2d}V_ie^{u_i} dx = \int_{\partial \gamma (B_{\epsilon}^+)} \partial_{\nu} u_i +o(\epsilon) + o(1)  \to \alpha_1 >0. $$

A contradiction.

\smallskip

Here we used a theorem of Hofmann see [11], which gives the fact that $ \gamma (B_{\epsilon}^+) $ is a Lipschitz domain. Also, we can see that $ \gamma ((-\epsilon, \epsilon)) $ and $ \gamma (\partial^+ B_{\epsilon}^+) $ are submanifolds.  

We start with a Lipschitz domain $ B_{\epsilon}^+ $ because it is convex and by the univalent and conformal map $ \gamma $ the image of this domain $ \gamma (B_{\epsilon}^+) $ is a Lipschitz domain and thus we can apply the integration by part and here we know the explicit formula of the unit outward normal it is the usual unit outward normal (normal to the tangent space of the boundary which we know explicitly because we have two submanifolds).

In the case of the disk $ D = \Omega $, it is sufficient to consider $ B(0,\epsilon)  \cap D $ which is a Lipschitz domain because it is convex (and not necessarily $ \gamma (B_{\epsilon}^+) $).

There is a version of the integration by part which is the Green-Riemann formula in dimension 2 on a domain $ \Omega $. This formula holds if we assume that there is a finite number of points $ y_1,..., y_m $ such that $ \partial \Omega - (y_1,..., y_m) $ is a $ C^1 $ manifold and for $ C^1 $ tests functions.

\smallskip

{\bf Remarks about the conformal map :}

1-It sufficient to prove  that $ \gamma_1((-\epsilon,\epsilon))=
\partial \Omega \cap \tilde \gamma_1(B_{\epsilon})= \partial \Omega \cap \tilde \gamma_1(B_{\epsilon})\cap \{ |abscissa|< \epsilon \} $, for $ \epsilon >0 $ small enough. Where $ \tilde \gamma_1 $ is the holomorphic extension of $ \gamma_1(t)=t+i\phi(t) $. For this, we argue by contradiction, we have for $ z_{\epsilon} \in B_{\epsilon} $, $\tilde \gamma_1(z_{\epsilon}) =(t_{\epsilon},\phi (t_{\epsilon})) $ for $ |t_{\epsilon}|\geq \epsilon $. Because $ \tilde \gamma_1 $ is injective on $ B_1 $ and $ \tilde \gamma_1 =\gamma_1=t+i\phi(t) $ on the real axis, we have necessirely $ |t_{\epsilon}|\geq 1 $. But, by continuity $ |\tilde \gamma_1(z_{\epsilon})|\to 0 $ because $ z_{\epsilon} \to 0 $. And, we use the fact that $ |\tilde \gamma_1(z_{\epsilon})|=|(t_{\epsilon},\phi(t_{\epsilon}))|\geq |t_{\epsilon}|\geq 1 $, to have a contradiction.) (This means  that for a small radius when the graph go out from the ball, it never retruns to the ball). (This fact imply that, when we have a curve which cut $\partial \Omega $ in $ \tilde \gamma_1(B_{\epsilon}) $ then the point have an abscissa such that $ |abscissa|< \epsilon $. This fact (by a contradiction with the fact $ \partial \bar \Omega = \partial \Omega $ and consider paths), imply that the image of the upper part of the ball is one side of the curve and the image of the lower part is in the other side of the curve.

2- Also, we can consider directly the coordinate $ T $ and change the function $ u(z)\to u(T) $ by $ z=T/a_1+x_0 $.

Set: $ \psi(\lambda_1,\lambda_2)\to M \in \Omega $ such that, $ \overrightarrow{x_0M}=\lambda_1 i_1'+\lambda_2j_1' $ where $ (i_1',j_1') $ is a basis such that $ i_1'=e^{-i\theta} i_1, j_1'=e^{-i\theta} j_1 $. And, $ \phi(x_1,x_2)\to M $ such that, $ \overrightarrow{OM}=x_1i_1+x_2i_2 $ the canonical basis $ (i_1,j_1) $. Then, we have two charts $ \phi $ and $ \psi $ and the complex affix $ T_M=\lambda_1+i\lambda_2 $ and $ z_M=x_1+ix_2 $ are such that  (transition map):

$$ T_M/e^{i\theta}+x_0=z_M=\phi^{-1} o \psi(\lambda_1,\lambda_2), $$

We have:

$$ \partial_{\lambda_1}=\cos \theta \partial_{x_1}+\sin \theta \partial_{x_2}, $$

$$ \partial_{\lambda_2}=-\sin \theta \partial_{x_1}+\cos \theta \partial_{x_2}, $$

Thus, the metric in the chart $\psi $ or coordinates $ (\lambda_1,\lambda_2) $ is : $ g_{ij}^{\lambda}=\delta_{ij} $ and the Laplacian in the two charts, $ \psi $ and $\phi $ are the usual Laplacian $ \partial_{\lambda_1\lambda_1}+\partial_{\lambda_2\lambda_2} $.

We write:

$$ \Delta u(M)=\Delta_{\lambda} (u o\psi(\lambda_1,\lambda_2)) $$

And then we apply the conformal map $\tilde \gamma_1 $ which send the affix $ T_M $, $ M $ in a neighborhood of $ x_0\in \partial \Omega $ to $ B_{\epsilon} $ with the fact that send $ T_M, M\in \partial \Omega $ to the real axis $ (-\epsilon,\epsilon) $ and the other parts of $ \Omega $ and $ \bar \Omega^c $.

3-We can remark that a definition of $ C^k, k \geq 1 $ domain, is equivalent to a definition of a submanifold with the condition $ \partial \bar \Omega = \partial \Omega $ or $ {\dot{\bar \Omega }} = \Omega $.

{\bf Remark 2: about a variational formulation.}

we consider a solutions in the sense of distribution. By the same argument (in the proof of the maximum principle obtained by Kato's inequality $ W^{1,1}_0 $ is sufficient), see [6], we prove that the solutions are in the sense $ C^2_0 $ of Agmon, see [1]. Also, we have corollary 1 of [7]. We use Agmon's regularity theorem. We return to the usual variational formulation in $ W^{1,2}_0 $ and thus we have the estimate of [8] or by Stampacchia duality theorem in $ W^{1,2}_0 $.

\end{document}